\documentclass[12 pt,a4paper]{article}
\usepackage{cite}
\usepackage{amsmath,amsfonts,amscd,amssymb}
\usepackage{extsizes}
\usepackage[linkcolor=blue, urlcolor=blue, citecolor=blue, colorlinks, bookmarks]{hyperref}
\usepackage[left=2cm,right=2cm,top=2cm,bottom=2cm]{geometry}
\usepackage{setspace}
\newtheorem{Theorem}{Theorem}[section]\numberwithin{Theorem}{section}
\newtheorem{Proposition}{Proposition}\numberwithin{Proposition}{section}

\newtheorem{Remark}{Remark}

\newtheorem{Lemma}[Theorem]{Lemma}
\numberwithin{equation}{section}

\begin{document}
	\title{\textbf{Cartan Connection for \textsl{h-}Matsumoto change}}
	\author{M.$\,$K.$\,$\textsc{Gupta},  Abha \textsc{Sahu}\thanks{Corresponding author}~,  Suman \textsc{Sharma}\\
	\normalsize{Department of Mathematics}\\[-3mm]
\normalsize{Guru Ghasidas Viswavidyalaya, Bilaspur (C.G.), India}\\[-3mm]
\small{E-mail: mkgiaps@gmail.com,  abhasahu118@gmail.com, sharma.suman209@gmail.com}}
	\date{}
	\maketitle
	\begin{abstract}
		In the present paper, we have studied the Matsumoto change $\overline{L}(x,y)= \frac{L^{2}(x,y)}{L(x,y) - \beta(x,y)} $ with an \textsl{h-}vector $b_{i}(x,y)$. We have derived some fundamental tensors for this transformation. We have also obtained the necessary and sufficient condition for which the Cartan connection coefficients for both the spaces $F^{n}=(M^n,L)$ and $\overline{F}^{\,n}=(M^{n},\overline{L})$ are same.
	\end{abstract}
	\section{Introduction}
		Let $M$ be an $ n $-dimensional $C^{\infty}$ Manifold and $T_{x}M$ denotes the tangent space of $M$ at $x$. The tangent bundle of $M$ is the union of tangent space $TM:=\underset{x\epsilon M}{\bigcup}T_{x}M$. A function $L: TM\rightarrow [0,\infty)$ is called Finsler metric function if it has the following properties\cite{shen2001lectures} \\[-1cm]
	\begin{enumerate}
		\item $L$ is $C^{\infty}$ on $TM \backslash\{0\}$,\\[-8mm]
		\item  For each $x\epsilon{M}$, $L_{x}:= L\vert_{T_{x}M}$ is a Minkowaski norm on $T_{x}M$.
	\end{enumerate}
	The pair $(M^{n}, L)$ is then called a Finsler space. The normalized supporting element, metric tensor, angular metric tensor and Cartan tensor are defined by $l_{i}=\dot{\partial_{i}}L$, $g_{ij}=\frac{1}{2}\dot{\partial_{i}}\dot{\partial_{j}}L^2 $, $ h_{ij}=L\dot{\partial_{i}}\dot{\partial_{j}}L $ and $ C_{ijk}=\frac{1}{2}\dot{\partial_{k}}g_{ij} $ respectively. The Cartan connection for the Finsler space $F^n$ is given by $(F^{i}_{jk},N^{i}_{j},C^{i}_{jk})$. The \textsl{h-}covariant and \textsl{v-}covariant derivative of the tensor $T^{\,i}_j$ with respect to Cartan connection, are respectively given as follows:\\[-6mm]
	\begin{equation*}\begin{split}
			T^{\,i}_{j|k}&=\delta_{k}T^{\,i}_{j} + T^{\,r}_{j}F^{i}_{rk} - T^{\,i}_{r}F^{r}_{jk}\,,\\
			T^{\,i}_{j}|_{k}&=\dot{\partial_{k}}T^{\,i}_{j} + T^{\,r}_{j}C^{i}_{rk} - T^{\,i}_{r}C^{r}_{jk}\,,
		\end{split}
	\end{equation*}
where $\delta_{k}$ is differential operator $\delta_{k}=\partial_{k}-N^{r}_{k}\dot{\partial_{r}}$. \\
	In 1984, C. Shibata \cite{shibata1984invariant} introduced the change $\overline{L}=f(L,\beta)$ as a generalization of Randers change, where $ f $ is positively homogeneous function of degree one in $ L $ and $ \beta(x,y)=b_{i}(x)y^{i} $. This change is called $ \beta $-change. An important class of $\beta$-change is Matsumoto change, given by 
	$\overline{L}(x,y) = \frac{L^2}{L-\beta}\,.$
If $L(x,y)$ reduces to a Riemannian metric then $\overline{L}(x,y)$ becomes Matsumoto metric. A famous example of Finsler space ``A slope measure of a moutain with respect to time measure'' was given by M.$\,$Matsumoto\cite{matsumoto1989slope}. Due to his great contribution in Finsler geometry, this metric was named after him.\\
	A.$\,$Tayebi et al. \cite{tayebi2014kropina} and Bankteswar Tiwari et al. \cite{tiwari2017generalized} discussed the Kropina change and generalized Kropina change respectively, for the Finsler space with $m^{th}$ root metric. In 2017, A.$\,$Tayebi et al.\cite{tayebi2017matsumoto} obtained the condition for the Finsler space given by Matsumoto change to be projectively related with the original Finsler space.\\[2mm]
	 The concept of \textsl{h-}vector $b_{i}$, was first introduced by H. Izumi \cite{izumi1980conformal}, which is \textsl{v-}covariant constant with respect to Cartan connection and satisfies 
		$L C^{h}_{ij} b_{h}= \rho h_{ij} \,,$
	 where $\rho$ is a non-zero scalar function. He showed that the scalar $\rho$ depends only on positional coordinates \textit{i.e.} $\dot{\partial_{i}}\rho =0$. From the definition of \textsl{h-}vector, it is clear that it depends not only on positional coordinates, but also on directional arguments.\\[2mm]
	Gupta and Pandey  \cite{gupta2009hypersurfaces , gupta2015finsler}, discussed certain properties of Randers change  and Kropina change with an \textsl{h-}vector. They\cite{gupta2015finsler} showed that 	\textit{If the \textsl{h-}vector is gradient then the scalar $\rho$ is constant}, \textit{i.e.} $\partial_{j}\rho=0$. In 2016, Gupta and Gupta \cite{gupta2017h, gupta2016hypersurface} have analysed Finsler space subjected to \textsl{h-}exponential change.\\
	 In the present paper, we have studied a Finsler metric defined by\\[-4mm]
	\begin{equation}\label{eq1}
		\overline{L}(x,y)= \frac{L^{2}(x,y)}{L(x,y) - b_{i}(x,y)y^{i}}\,,
	\end{equation}
	where $b_{i}(x,y)$ is an \textsl{h}-vector in $(M^{n},L)$.\\
	The structure of this paper is as follows: In section $2$, we have obtained the expressions for different fundamental tensors of the transformed Finsler space. In section $3$, we have observed how the Cartan connection coefficients change due to Matsumoto change  with an \textsl{h}-vector and also find the necessary and sufficient condition for which both connection coefficients would be same. 
	\begin{Remark}
		H.$\,$S.$\,$Shukla et al. \emph{\cite{shuklamatsumoto}} also discussed Matsumoto change of Finsler metric by \textsl{h}-vector. Unfortunately, the results are wrong because of wrong computation in Lemma $1.1$ of\emph{\cite{shuklamatsumoto}}.
	\end{Remark}
	\section{The Finsler space $\overline{F}^{\,n}= (M^{n},\overline{L})$}
	Let the Finsler space transformed by the Matsumoto change \eqref{eq1} with an \textsl{h}-vector, be denoted by $\overline{F}^{\,n}= (M^{n},\overline{L})$. If we denote $\beta=b_{i}(x,y)y^{i}$, then indicatory property of  angular metric tensor yields $\dot{\partial_{j}}\beta=b_{j}\,$. Throughout this paper, we have barred the geometrical objects associated with $\overline{F}^{\,n}$.
	\\From \eqref{eq1}, we get the normalized supporting element as \\[-6mm]
	\begin{equation}\label{eq2}
		\overline{l}_{i}= \frac{\tau}{(\tau-1)}\,l_{i}+ \frac{\tau^{2}}{(\tau-1)^{2}}\,m_{i}\,,
	\end{equation} \\[-1cm]
where $\tau= \frac{L}{\beta}$ and $m_{i}= b_{i} - \frac{1}{\tau}l_{i}\,.$
\begin{Remark}
The covariant vector  $m_{i}$ statisfies the following relations  \\
\emph{(i)}$\:\: m_{i}\neq 0 \qquad$\emph{(ii)}$\:\: m^{i}=g^{ij}m_{j}\qquad$ \emph{(iii)}$\:\: m^{i}m_{i}=b^2-\frac{1}{\tau^2}=m^2\qquad$\emph{(iv)}$ \:\: m_{i}y^{i}=0 \,.$	
\end{Remark}
Differentiating equation \eqref{eq2} with respect to $y^{j}$, and using the notation $L_{ij}=\dot{\partial_{j}}l_{i}$ we get 
	\begin{equation*}
		\overline{L}_{ij}= \frac{\tau(\tau+\rho\tau-2)}{(\tau-1)^{2}}\, L_{ij}+\frac{2\tau^{2}}{\beta(\tau-1)^3}\,m_{i}m_{j}.
	\end{equation*}
Therefore, the angular metric tensor $ \overline{h}_{ij}$ is obtained as
\begin{equation}
	\overline{h}_{ij}=\frac{\tau^{2}(\tau+\rho\tau-2)}{(\tau-1)^{3}}\,h_{ij} + \frac{2\tau^4}{(\tau-1)^4}m_{i}m_{j}\,.
\end{equation}
	The metric tensor $\overline{g}^{}_{ij}=\overline{h}_{ij}+\overline{l}_{i}\overline{l}_{j}$ is given by
	\begin{equation}
		\overline{g}^{}_{ij}=\frac{\tau^{2}(\tau+\rho\tau-2)}{(\tau-1)^{3}}\,g^{}_{ij}+\frac{\tau^{2}(1-\rho\tau)}{(\tau-1)^{3}}\,l_{i}l_{j}+\frac{\tau^{3}}{(\tau-1)^{3}}\,(m_{i}l_{j}+m_{j}l_{i})+ \frac{3\tau^4}{(\tau-1)^4}\,m_{i}m_{j}\,,
	\end{equation}
which can be rewritten as\\[-5mm]
	\begin{equation}\label{eq3}
		\overline{g}^{}_{ij}= p\,g^{}_{ij}+ p^{}_{1}l_{i}l_{j}+ p^{}_{2}(m_{i}l_{j}+m_{j}l_{i})+p^{}_{3}\,m_{i}m_{j}\,,
	\end{equation}
 where
 \begin{equation*}
		p=\frac{\tau^{2}(\tau+\rho\tau-2)}{(\tau-1)^{3}}, \quad p^{}_{1}=\frac{\tau^{2}(1-\rho\tau)}{(\tau-1)^{3}},
	\quad	 p^{}_{2}=\frac{\tau^{3}}{(\tau-1)^{3}}, \quad p^{}_{3}= \frac{3\tau^4}{(\tau-1)^4}\,.
	\end{equation*}
\noindent The following lemma helps us to compute the inverse of metric tensor $\overline{g}^{}_{ij}\,$.
\begin{Lemma}
\emph{\cite{matsumoto1972c}:} Let $(m_{ij})$ be a non-singular matrix and $l^{}_{ij}= m_{ij}+n_{i}n_{j}$. The elements $l^{ij}$ of the inverse matrix, and the determinant of the matrix $(l^{}_{ij})$ are given by \\[-6mm]
	\begin{equation*}
		l^{ij}= m^{ij}-(1+n_{k}n^{k})^{-1}n^{i}n^{j},\quad det(l^{}_{ij})=(1+n_{k}n^{k})det(m_{ij})
	\end{equation*}
respectively, where $m^{ij}$ are elements of the inverse matrix of $(m_{ij})$ and $n^{k}=m^{ki} n_{i}$.
\end{Lemma}
 The inverse metric tensor of $\overline{F}^{\,n}$ can be derived as follows:\\[-6mm]
 \begin{equation}\label{eq4}
		\overline{g}^{\,ij}= q\,g^{ij}+ q_{1}\,l^{i}l^{j}+ q_{2}\,(l^{i}m^{j}+m^{i}l^{j})+ q_{3}\,m^{i}m^{j}\,,
	\end{equation}\\[-12mm]
	where\\[-6mm]
	\begin{equation*}
		  \quad q= \frac{1}{p}\,,\qquad q_{1}= \frac{-1}{2}\Big[\frac{p_{1}{p^{}_{3}}-p^{2}_{2}}{(p_{1}+p)p_{3}-p_{2}^{2}}+\frac{2p^{2}p^{2}_{2}p_{3}}{(3p+2p_{3}m^{2})\{(p_{1}+p)p_{3}-p_{2}^{2}\}^{2}}\Big] ,
		\end{equation*}	
\begin{equation*}
			 q_{2}=\frac{-2p_{2}p_{3}}{(3p+2p_{3}m^{2})\{(p_{1}+p)p_{3}-p_{2}^{2}\}}, \qquad q_{3}=\frac{-2p_{3}}{p(3p+2p_{3}m^{2})}\,.
\end{equation*}	
	The Cartan tensor $\overline{C}_{ijk}$ is obtained by differentiating the equation (\ref{eq3}) with respect to $ y^{k} $, as follows:\\[-8mm]
	\begin{equation}\label{eq6}
		\overline{C}^{}_{ijk}=p\,C^{}_{ijk}+V^{}_{ijk}\,,
	\end{equation}\\[-1cm]
where\\[-12mm]
	\begin{equation*}
	V_{ijk}=K_{1}(h_{ij}m_{k}+h_{jk}m_{i}+h_{ki}m_{j})+K_{2}\,m_{i}m_{j}m_{k}
	\end{equation*}
and\\[-12mm]
	\begin{equation*}
		K_{1}= \frac{\tau^{3}(\tau+3\rho\tau-4)}{2L(\tau-1)^4},\hspace{.3cm} K_{2}= \frac{6\tau^{4}}{\beta(\tau-1)^5}\,.
	\end{equation*}\\[-16mm]
\begin{Remark}
 From above we can retrieve  relations between the scalars as
\begin{equation*}
	\frac{\partial\,p}{\partial\, \tau}=-\frac{2L}{{\tau}^2}\,K_{1}\,, \quad 	\frac{\partial\,p_{3}}{\partial\, \tau}=-\frac{2L}{{\tau}^2}\,K_{2}\,,	
\end{equation*}
\begin{equation*}
	K_{1}=\frac{1}{2L}\left\lbrace p_{2} + p_{\,3}\left(\rho - \frac{1}{\tau} \right)  \right\rbrace\,\quad \emph{and} \quad p_{\,1} + p_{2}\left(\rho - \frac{1}{\tau} \right)  = 0\,.
\end{equation*}
\end{Remark}
	From equation \eqref{eq4} and \eqref{eq6}, we get the \textsl{(h)hv-}torsion tensor $\overline{C}^{\,i}_{jk}$\\[-6mm]
	\begin{equation}\label{eq7}
		\overline{C}^{\,i}_{jk}= C^{\,i}_{jk}+M^{\,i}_{jk}\,,
	\end{equation}where\\[-12mm]
	\begin{equation*}\begin{split}
			 \quad M^{i}_{jk}&= q\,K_{1}(m_{k}h^{i}_{j}+m_{j}h^{i}_{k})+(q_{2}\,l^{i}+q_{3}\,m^{i})\left\lbrace 2K_{1}m_{j}m_{k}+\frac{p}{L} \rho \,h^{}_{jk}\right\rbrace \\
			&+\left\lbrace q\,m^{i}+(q_{2}\,l^{i}+q_{3}\,m^{i})m^{2}\right\rbrace \left( K_{2}m_{j}m_{k}+K_{1}h^{}_{jk}\right).
		\end{split}
	\end{equation*}
	\section{Cartan Connection of the space $\overline{F}^{\,n}$}
	The Cartan connection for a Finsler space $\overline{F}^{\,n}$ is given by the  traid $(\overline{F}^{\,i}_{jk},\overline{N}^{\,i}_{j},\overline{C}^{\,i}_{jk})$.
	The \textsl{v-}connection coefficient $ \overline{C}^{\,i}_{jk} $ is given by equation \eqref{eq7}. Now, we are obtaining the \textsl{h-}connection coeffiecient $ \overline{F}^{\,i}_{jk} $ and non-linear connection coeffiecient $ \overline{N}^{\,i}_{j} $.\\ First, we will try to find canonical spray of the transformed space $ \overline{F}^{\,n} $.\\
Differentiating equation \eqref{eq3} with respect to $x^{k}$, and using the definition of \textsl{h-}covariant derivative, we obtain\\
\begin{equation}\label{eq8}
	\begin{split}
		\partial^{}_{k} \overline{g}^{}_{ij} = &\,p\,\partial^{}_{k} g^{}_{ij}+p^{}_{1}(l_{i}l_{r}F^{r}_{jk}+l_{j}l_{r}F^{r}_{ik})
		+p^{}_{2}(\rho^{}_{k}h^{}_{ij}+l^{}_{i}b^{}_{j|k}+l^{}_{j}b^{}_{i|k}+m^{}_{r}F^{r}_{jk}l_{i}+m^{}_{r}F^{r}_{ik}l_{j}\\&+m^{}_{i}F^{r}_{jk}l_{r}+m^{}_{j}F^{r}_{ik}l_{r})+p^{}_{3}(m_{i}b^{}_{j|k}+m_{j}b^{}_{i|k}+m_{i}m_{r}F^{r}_{jk}+m_{j}m_{r}F^{r}_{ik})\\&+ 2(K_{1}h^{}_{ij}+2K_{2}m_{i}m_{j})(\beta^{}_{k}+ N^{r}_{k}m_{r})+K_{1}(h_{jr}N^{r}_{k}m_{i}+h_{ir}N^{r}_{k}m_{j})\,,
	\end{split}
\end{equation} where $\partial_{k} \rho={\rho}_{|k}=\rho^{}_{k}$ and ${\beta}_{|k}={\beta}_{k}$.\\
\noindent Applying Christoffel process with respect to indices $i\,,j\,,k$ in above equation, we obtain the coefficient of Christoffel symbol as follows:\\[-4mm]
\begin{equation}\begin{split}
		\overline{\gamma}_{ijk} = p{\gamma}_{ijk}&+ \mathfrak{S}_{ijk}\left\lbrace \frac{p^{}_{2}}{2}\rho_{k}h_{ij}+\left( \beta_{k} + N^{r}_{k}m_{r}\right) B_{ij}\right\rbrace  + Q_{i}F_{jk}+Q_{k}F_{ji} +Q_{j}E_{ik}\\&+(\overline{g}_{rj}-p\,g_{rj})\Big\{{\gamma}^{\,r}_{ik} + g^{rt}(C_{ikm}N^{m}_{t}-C_{tkm}N^{m}_{i}-C_{itm}N^{m}_{k})\Big\}\,,
\end{split}\end{equation}\vspace{2mm}
where the symbol $\mathfrak{S}_{ijk}$ is defined as $\mathfrak{S}_{ijk}\,U_{ijk}=U_{ijk}-U_{jki}+U_{kij}$
and we have used the notation\\[-6mm]
\begin{equation*}
Q_{i}=p^{}_{2}l_{i}+p^{}_{3}m_{i}\,, \quad B_{ij}=K^{}_{1}h_{ij}+K^{}_{2}m_{i}m_{j}\,,	
\end{equation*}
\begin{equation*}
		2E_{ij}=b_{i|j}+b_{j|i}\,, \!\qquad 2F_{ij}=b_{i|j}-b_{j|i}\,.
\end{equation*}	\\[-2cm]
\begin{Remark}
	The tensors $Q_{i}$ and $ B_{ij} $ statisfy the following\\[2mm]
	\emph{(i)} $\quad \dot{\partial_{j}}Q_{i}=B_{ij}\qquad$
	\emph{(ii)}$\quad B_{ij}= B_{ji}\qquad$
	\emph{(iii)}$ \quad B_{ij}y^{i}=0 \,.$
\end{Remark}	
The Christoffel Symbol of second kind of the Finsler space $\overline{F}^{\,n}$ is given by\\[-4mm]
\begin{equation}\begin{split}\label{eq9}
		\overline{\gamma}^{\,i}_{jk}={\gamma}^{\,i}_{jk} &+ (g^{it}-p\overline{g}^{\,it})(C_{jkm}N^{m}_{t} - C_{tkm}N^{m}_{j} - C_{jtm}N^{m}_{k}) \\[2mm]&+ \overline{g}^{\,is}\mathfrak{S}_{jsk}\left[ \left\lbrace \frac{p^{}_{2}}{2}\rho_{k}h_{js} + (\beta_{k} + N^{r}_{k}m_{r})B_{sj}\right\rbrace  + Q_{j}F_{sk}  +  Q_{k}F_{sj}+  Q_{s}E_{jk} \right] \,.
\end{split}\end{equation}\\[-3mm]
Transvecting equation \eqref{eq9} by $ y^{j}y^{k} $ and using ${G}^{\,i}= \frac{1}{2}{\gamma}^{\,i}_{jk}y^{j}y^{k}$, we get\\[-6mm]
\begin{equation}\label{eq15}
	\overline{G}^{\,i}={G}^{\,i} + {D}^{\,i}\,,
\end{equation}\vspace{-4mm}
where  
\begin{equation}\label{eq20}
{D}^{\,i}=\frac{1}{2}\,\overline{g}^{\,is}\big[Q_{s}E_{oo} + 2p_2LF_{so}\big]\,.
\end{equation}
Thus, we have:
\begin{Proposition}
	The spray coefficient of the transformed space is given by equation \eqref{eq15}.
\end{Proposition}
\begin{Remark}
	In the subscript zero `o' is used to denote the transvection by $y^{i}$, \textit{ i.e.} $ F_{so}=F_{si}y^{i} $.	
\end{Remark}
Differentiating  equation \eqref{eq15} with respect to $y^{i}$ and using $\dot{\partial_{j}}G^{i}=N^{i}_{j}$ and $\dot{\partial_{j}}\overline{g}^{\,is}= -2\,\overline{g}^{\,ir}\overline{C}^{\,s}_{rj}$, we get 
\begin{equation}\label{eq12}
	\overline{N}^{\,i}_{j}={N}^{\,i}_{j} + {D}^{\,i}_{j}\,,
\end{equation}\\[-12mm]
	where \\[-8mm]
\begin{equation}\label{eq30}
	 {D}^{\,i}_{j}= \overline{g}^{\,ir} \Big\{ -2D^{m}(\,p\,C_{mrj} + V_{mrj}) + Q_{r}E_{oj} + E_{oo}B_{rj} +p_2LF_{rj} + Q_{j}F_{ro} + \frac{p_2}{2}\rho^{}_{k}y^{k}h_{rj}\Big\} .
\end{equation}
Thus, we have:
\begin{Proposition}
	The non-linear connection coefficient of the transformed space is given by the equation \eqref{eq12}.
\end{Proposition}
Now, we are in a position to obtain the Cartan connection coefficient for the space  $\overline{F}^{\,n}$. We know that the relation between the Christoffel symbol and Cartan connection coefficient is given by
\begin{equation*}
	{F}^{\,i}_{jk} = {\gamma}^{\,i}_{jk} + {g}^{\,is}({C}_{jkr}{N}^{r}_{s} - {C}_{skr}{N}^{r}_{j} - {C}_{jsr}{N}^{r}_{k})\,.
\end{equation*}
In view of equation \eqref{eq6}, \eqref{eq9} and \eqref{eq12}, we have
\begin{equation*}\begin{split}
		\overline{F}^{\,i}_{jk} = {\gamma}^{\,i}_{jk} &+ (g^{it}-p\overline{g}^{\,it})(C_{ikm}N^{m}_{t} - C_{tkm}N^{m}_{i} - C_{itm}N^{m}_{k}) \\ &+ \overline{g}^{\,is}\left\lbrace \mathfrak{S}_{jsk}\left(  {\frac{p_2}{2}\rho^{}_{k}h_{js} + (\beta_{k} + N^{r}_{k}m_{r})B_{sj}}\right)  + Q_{j}F_{sk} +Q_{k}F_{sj}+  Q_{s}E_{jk}   \right\rbrace \\ &+ \overline{g}^{\,is}\Big\{(pC_{jkr} + V_{jkr})(N^{r}_{s} + D^{r}_{s}) - (pC_{skr} + V_{skr})(N^{r}_{j} + D^{r}_{j}) - (pC_{jsr} + V_{jsr})(N^{r}_{k} + D^{r}_{k})\Big\}
\end{split}\end{equation*}
which can be simplifed as\\[-5mm]
\begin{equation*}\begin{split}
		\overline{F}^{\,i}_{jk} = {F}^{\,i}_{jk} &+ \overline{g}^{\,is}\Big\{ \mathfrak{S}_{jsk}\Big(\frac{p_2}{2}\rho^{}_{k}h_{js} + \beta_{k}B_{js} - p\, C_{jsr}D^{r}_{k} - V_{jsr}D^{r}_{k}\Big) + Q_{j}F_{sk} + Q_{k}F_{js} +Q_{s}E_{jk} \Big\}.
\end{split}\end{equation*}
Above equation can be rewritten as\\[-8mm]
\begin{equation}\label{eq16}
	\overline{F}^{\,i}_{jk} = {F}^{\,i}_{jk} + {D}^{\,i}_{jk}\,,
\end{equation}\\[-12mm]
	where\\[-8mm]
\begin{equation}\begin{split}\label{eq10}
	 {D}^{\,i}_{jk}=\overline{g}^{\,is}\Big\{ \mathfrak{S}_{jsk}\Big(\frac{p_2}{2}\rho^{}_{k}h_{js} + \beta_{k}B_{js} - p\, C_{jsr}D^{r}_{k} - V_{jsr}D^{r}_{k}\Big) + Q_{j}F_{sk} + Q_{k}F_{sj} + Q_{s}E_{jk}  \Big\}\,.
\end{split}\end{equation}
Hence, we have:
\begin{Theorem}
The relation between the Cartan connection coefficeint of $F^{n}$ and  $\overline{F}^{\,n}$ is given by  equation \eqref{eq16}.
\end{Theorem}
\begin{Remark}\label{R1}
	The tensors $D^{i}_{jk},\: D^{i}_{j}$ and $D^{i}$ are related as\\
	\emph{(i)}$  \quad D^{\,i}_{jk}\,y^{k}=D^{\,i}_{j}\,, \qquad$
\emph{(ii)}$ \quad D^{\,i}_{j}\,y^{j}=2D^{\,i}\,, \qquad  $
\emph{(iii)}$ \quad \dot{\partial_{j}}D^{\,i}=D^{\,i}_{j}.$
\end{Remark}
Now, we want to find the condition for which the Cartan connection coefficients for both spaces $F^{n}$ and $\overline{F}^{\,n}$ are same, \textit{i.e.} $\overline{F}^{\,i}_{jk} = {F}^{\,i}_{jk}$ then ${D}^{\,i}_{jk}=0$, which implies ${D}^{\,i}_{j}=0$, then ${D}^{\,i}=0$.
Therefore the equation \eqref{eq20} gives \\[-8mm]
\begin{equation*}
	2\,p_2LF_{io} + E_{oo}Q_{i}=0\,,
\end{equation*} \\[-8mm]
which on transvection by $y^{i}$ gives $E_{oo} = 0$ and then $F_{io}=0$.
Diffrentiating $E_{oo} = 0$ partially with respect to $y^{i}$ gives $E_{io}=0$. Therefore we have $E_{io}=0=F_{io}$, which implies $ b_{i|o}=b_{o|i}=\beta_{\,|i}=0$. Differentiating $ \beta_{\,|i} $  partially with respect to $ y^{j} $ and using the commutation formula $ \dot{\partial_{j}}(\beta_{\,|i})-(\dot{\partial_{j}}\beta)_{|i}=(\dot{\partial_{r}}\beta)C^{r}_{ij|o}\, $, we get \\[-6mm]
\begin{equation}\label{eq31}
	b_{j|i}=-\,b_{r}C^{r}_{ij|o}\,.
\end{equation} This will give us $F_{ij}\!=\!0$. Taking \textsl{h-}covariant derivative of $LC^{r}_{ij}b_{r}=\rho h_{ij}$ and using ${\rho}^{}_{|k}\!=\!0$, $L_{|k}=0$ and $h_{ij|k}=0$, we get \\[-6mm]
\begin{equation*}
\left( C^{r}_{ij}b_{r}\right)_{|k}=\left( \frac{\rho}{L} h_{ij}\right)_{|k}=0 \,.
\end{equation*}\\[-1cm]
This gives \\[-12mm]
\begin{equation*}
C^{r}_{sj}b_{r|k} + C^{r}_{sj|k}b_{r}=0\,.
\end{equation*}
Transvecting by $ y^{k}$ and using $b_{r|o}=0$, we get $ C^{r}_{ij|o}b_{r} =0$ and then equation \eqref{eq31} gives $b_{i|j}=0 $, \textit{i.e.} the \textsl{h-}vector $b_{i}$ is parallel with respect to Cartan connection of $F^{n}$.\\[3mm]
\underline{Conversely}, If $b_{i|j}=0$ then we get $E_{ij}=F_{ij}=0$ and $\beta_{i}=\beta_{\,|i}=b_{j|i}\,y^{j}=0$. Then equation \eqref{eq20} reduces to $D^{i}=0$. From $F_{ij}=0$ we have $ {\rho_{i}}=0 $, which implies $ D^{i}_{j}=0$.
 Therefore, from equation \eqref{eq10}, we get $D^{i}_{jk}=0$, which gives $\overline{F}^{\,i}_{jk} = {F}^{\,i}_{jk}$. Thus, we have:
\begin{Theorem}\label{T1}
	For the Matsumoto change with an \textsl{h-}vector, the Cartan connection coefficients for both spaces $F^{n}$ and $\overline{F}^{\,n}$ are the same if and only if the \textsl{h-}vector $b_{i}$ is parallel with respect to the Cartan connection of $F^{n}$.
\end{Theorem}
Now, differentiating equation \eqref{eq12} with respect to $ y^{k}$, and using $\dot{\partial_{k}}N^{i}_{j}=G^{i}_{jk}$, we obtain\\[-6mm]
\begin{equation}\label{eq26} 
	\overline{G}^{\,i}_{jk}=G^{i}_{jk}+\dot{\partial_{k}}\,D^{i}_{j}\,,
\end{equation}
where $ G^{i}_{jk} $ are the Berwald connection coeffiecients.\\
Now, if the \textsl{h-}vector $b_{i}$ is parallel with respect to the Cartan connection of $F^{n}$, then by the Theorem \ref{T1}, the cartan connection coefficients for both Finsler space $ F^{n} $ and $ \overline{F}^{\,n} $ are the same, \textit{i.e.} $D^{i}_{jk}=0$ which implies $D^{i}_{j}=0$. Then from equation \eqref{eq26}, we get $\overline{G}^{\,i}_{jk}=G^{i}_{jk}\,.$\\[2mm]
\underline{Conversely}$\,$, If $\overline{G}^{\,i}_{jk}=G^{i}_{jk}\,$ then, from equation \eqref{eq26}, we have $\dot{\partial_{k}}\,D^{i}_{j}=0$, which on transvecting by $ y^{j}$ and using Remark \ref{R1}, gives $D^{i}_{k}=0\,$. Using the same procedure as in the Theorem \ref{T1}, we get $b_{i|j}=0 $, \textit{i.e.} the \textsl{h-}vector $b_{i}$ is parallel with respect to Cartan connection of $F^{n}$. \\
Thus, we have:
\begin{Theorem}
	For the Matsumoto change with an \textsl{h-}vector, the Berwald connection coefficients for both spaces $F^{n}$ and $\overline{F}^{\,n}$ are the same if and only if the \textsl{h-}vector $b_{i}$ is parallel with respect to the Cartan connection of $F^{n}$.
\end{Theorem}
\section*{Conclusion}
In the present paper, The Cartan connection of the changed Finlser space is discovered and with the condition (\textsl{h-}vector $b_{i}$ is parallel, \textit{i.e.} $ b_{i|j}=0 \,$), the Cartan connection of both the spaces are same.\\
\textit{ For this transformation we can also find some geometric properties for the transformed Finsler space like the curvature tensor, torsion tensor, T-tensor etc.}\\
 Gupta and Pandey \cite{gupta2015finsler} have proved that, 
 	``$\,$For the Kropina change with an \textsl{h-}vector, the Cartan connection coefficients for both spaces $F^{n}$ and $\overline{F}^{\,n}$ are the same if and only if the \textsl{h-}vector $b_{i}$ is parallel with respect to the Cartan connection of $F^{n}\,$".
 We here observe that the Kropina change has finite number of terms whereas, Matsumoto change has infinite number of terms, although in both cases (finite and infinite) same result holds.\\
 \textit{\bfseries The goal for future study in this area is to identify a class of change with an \textsl{h-}vector  $b_{i}$ is parallel, for which the Cartan connection of both the Finsler space are same}. 
\normalsize

\end{document}